\documentstyle[11pt]{article}
\input amssymb.sty
\newcommand{\G}{\Gamma}
\newcommand{\p}{\partial}

\newcommand{\ti}[1]{\tilde{#1}}

\newcommand{\Om}{\Omega}
\newcommand{\de}{\delta}
\newcommand{\al}{\alpha}

\newcommand{\Te}{\Theta}

\newcommand{\D}{\Delta}
\newcommand{\ve}{\varepsilon}

\newcommand{\ze}{\zeta}
\newcommand{\si}{\sigma}

\newcommand{\clK}{{\cal K}}
\newcommand{\bfE}{{\bf E}}
\newcommand{\bfe}{{\bf e}}
\newcommand{\clZ}{{\cal Z}}

\newcommand{\clP}{{\cal P}}
\newcommand{\clF}{{\cal F}}
\newcommand{\clW}{{\cal W}}

\newcommand{\mat}[4]{\left(\begin{array}{cc}{#1}&{#2}\\{#3}&{#4}
\end{array}\right)}
\newcommand{\oh}{\frac{1}{2}}

\newcommand{\SLR}{{\rm SL}(2,{\mathbb R})}

\newcommand{\SL}{{\rm SL}(2,{\mathbb C})}

\def\f1#1{\frac{1}{#1}}

\newcommand{\bz}{\bar{z}}

\newcommand{\bxi}{\bar{\xi}}
\def\mC{{\mathbb C}}
\def\mZ{{\mathbb Z}}

\def\mN{{\mathbb N}}

\newcommand{\ot}{\otimes}
\newtheorem{predl}{Proposition}[section]
\newtheorem{defi}{Definition}[section]

\renewcommand{\theequation}{\thesection.\arabic{equation}}
\newcommand{\beq}[1]{\begin{equation}\label{#1}}
\newcommand{\eq}{\end{equation}}

\def\mapright#1{\smash{
\mathop{\longrightarrow}\limits^{#1}}}
\def\mapleft#1{\smash{
\mathop{\longleftarrow}\limits^{#1}}}

\def\mapdown#1{\Big\downarrow\rlap
{$\vcenter{\hbox{$\scriptstyle#1$}}$}}

\begin{document}
\vspace{10mm}
\vspace{0.3in}
\begin{flushright}
 ITEP-TH-17/03\\
\today
\end{flushright}
\vspace{10mm}
\begin{center}
{\Large \bf Poisson formula for a family of non-commutative Lobachevsky spaces}\\
\vspace{5mm}
M.A.Olshanetsky\\
{\sf ITEP, 117259, Moscow, Russia}\\
{\em olshanet@gate.itep.ru}\\
V.-B.K.Rogov \footnote{The work of the second author was supported by
the NIOKR MPS RF Foundation}\\
{\sf MIIT, 127994, Moscow, Russia} \\
{\em vrogov@cemi.rssi.ru}\\
\vspace{5mm}
\end{center}
\begin{flushright}
{\sl Dedicated to the memory of Fridrikh Izrailevich Karpelevich }
\end{flushright}
\vspace{10mm}
\begin{abstract}
We define an analog of the Poisson integral formula
for a family of the non-commutative Lobachevsky spaces. The $q$-Fourier
transform of the Poisson kernel is expressed through the
$q$-Bessel-Macdonald function.
\end{abstract}

\section{Introduction}
\setcounter{equation}{0}

The classical Lobachevsky space ${\bf L}^3$ can be identified with 3d hyperboloid
\beq{em}
{\bf L}^3=\{x_0^2-x_1^2-x_2^2-x_3^2=1,~x_0>0\}\,,
\eq
equipped with the hyperbolic metric. It can be represented as the set of
the second order positive definite Hermitian
matrices
$$
x=
\left(
\begin{array}{cc}
x_0-x_1 & x_2-ix_3 \\
x_2+ix_3 & x_0+x_1
\end{array}\right)\,,~~\det x=1\,.
$$
It means that ${\bf L}^3$ is the quotient space
$SU_2\backslash SL_2({\mC})$. In this way
the classical Lobachevsky space is a particular example of
the symmetric spaces.

We introduce the horospheric coordinates $(H,z,\bz)$ on ${\bf L}^3$
$$
\left\{
\begin{array}{l}
x_0=\oh (|z|^2H+H+H^{-1})\,,\\x_1=\oh(- |z|^2H+H-H^{-1})\,,\\
x_2=\oh H(z+\bz)\,,\\ x_3=-\frac{i}{2}H(z-\bz)\,.
\end{array}
\right.
$$

Let $\Om$ be the Laplace-Beltrami operator on ${\bf L}^3$ corresponding to the metric
that comes from the embedding (\ref{em}). In the horospheric coordinates it
has the form
\beq{Om}
\Omega=\frac12H^2\frac{\p^2}{\p H^2}+\frac32H\frac\p{\p H}+
2H^{-2}\frac{\p^2}{\p{\bar z}\p z}\,.
\eq
Consider the equation
\beq{3.4}
\left(\frac12\Omega+\frac14\right)F_\nu(\bar z,H,z)=
\frac{\nu^2}4F_\nu(\bar z,H,z),\quad \nu\geq 0.
\eq

The non-negative solutions of this equation are described in the following
way.
 Let
$$
{\bf C}^3=\{x_0^2-x_1^2-x_2^2-x_3^2=0,~x_0>0\}\,.
$$
be the upper pole of the cone. It can be identified with the Hermitian matrices with
$\det x=0$ and $x_0\geq 0$.
The horospheric coordinates  $(\al,\xi,\bxi)$ on ${\bf C}^3$ take the form
$$
\left\{
\begin{array}{l}
x_0=\oh\al(|\xi|^2+1)\,,\\x_1=\oh\al(1-|\xi|^2)\,,\\
x_2=\frac{\al}{2}(\xi+\bxi)\,,\\x_3=-\frac{i\al}{2}(\xi-\bxi)\,.
\end{array}
\right.
$$
We identify {\sl the absolute} $\Xi$ of ${\bf L}^3$
with the section
$\al=const$ of ${\bf C}^3$ completed with the point $(z=\infty)$.
Then  solutions of (\ref{3.4}) can be represented as the Poisson integral
\beq{pf}
F_\nu(\bar z,H,z)=\int_\Xi{\cal P}(\bxi-\bz,H,\xi-z)d\mu(\xi,\bxi)\,,
\eq
where
\beq{pk}
{\cal P}(\bxi-\bz,H,\xi-z)= (\{(H^{-1}+|z-\xi|^2H)^{-\nu-1}
\eq
is {\sl the Poisson kernel} and $d\mu(\xi,\bxi)=\phi(\xi,\bxi)d\xi\bxi$ is a
Borel distribution $\Xi$.
The classical theorem states that:\\
i. non-negative solutions
exist for $\nu\geq 0$ and have the representation (\ref{pf});\\
ii. $\phi(\xi,\bxi)$ is
uniquely reconstructed from $ F_\nu(\bar z,H,z)$.

This theorem was generalized
by Karpelevich and Furstenberg on  arbitrary symmetric spaces of the non-compact
type \cite{Ka,Fu}.

We investigate an analog of this representation for a family of non-commutative
Lobachevsky spaces constructed in Ref.\,\cite{OR3}.
We use also another form of (\ref{pf}).
Let $\Phi_\nu(\bar s,H,s)$ be the
Fourier transform of $F_\nu(\bar z,H,z)$
with respect to $(\bz,z)$. It satisfies the ordinary differential equation
\beq{3.5}
\left(\frac14H^2\frac{d^2}{dH^2}+\frac34H\frac{d}{dH}-
H^{-2}\bar ss+\frac14\right)\Phi_\nu(\bar s,H,s)=
\frac{\nu^2}4\Phi_\nu(\bar s,H,s)\,.
\eq
The solutions to (\ref{3.5}) decreasing for $H\to 0$ are the
functions
\beq{3.6}
\Phi_\nu(\bar s,H,s)=\frac\pi{\G(\nu+1)}H^{-1}
K_\nu(2\sqrt{\bar ss}H^{-1})(\bar ss)^{\frac\nu2}\psi(\bar s,s),
\eq
where $K_\nu$ is the Bessel-Macdonald function, and
$\psi_\nu(\bar s,s)$ is the Fourier transform of $\phi(\xi,\bxi)$.

The family of the non-commutative Lobachevsky spaces depends on the
deformation parameter $0<q<1$ and $\de=0,1,2$. The classical (commutative) limit
corresponds to $q=1$. The discrete parameter $\de$ is responsible for
the form of the Laplace-Beltrami operator.
Our results only partly reproduce the classical situation.
We  just construct an analog of the both representations
(\ref{pf}) and (\ref{3.6}) for a
space of functions on the non-commutative absolute and its Fourier dual.
In the limit $q\to 1$ we come to (\ref{pf}) and (\ref{3.6}).
Some kernels on the unit disc
were considered in Ref.\,\cite{Va}

\section{ Non-commutative Lobachevsky spaces (NLS). }
\setcounter{equation}{0}

{\bf 1. General definition.}

The  description of  NLS is based on an analog of the horospheric coordinates.
 Let ${\bf L}_{\de,q}$ be an associative $*$-algebra over
${\mathbb C}$ with the unit and three generators
$$
(z^*,H,z)\,,~~H^*=H\,,~(z)^*=z^*\,,
$$
and the commutation relations depending on two parameters $q\in (0,1)$ and
$\de=0, 1, 2$
\beq{4.4}
Hz=q^\de zH\,,~~~~z^*H=q^\de Hz^*\,,
~~~~z^*z=q^{2-2\de}zz^*-q^{-\de}(1-q^2)H^{-2}\,.
\eq
To eliminate the ambiguities related to the non-commutativity
we consider only the ordered monomials putting $z^*$ on the left side, $z$ on
the right side and keeping $H$ in the middle of the monomials:
$$
w(m,k,n)=(z^*)^{m}H^kz^n\,.
$$
The symbol $\ddag f(z^*,H,z)\ddag$ denotes that all monomials are
ordered
$$
\ddag f(z^*,H,z)\ddag=\sum_{m,k,n=-\infty}^\infty a_{m,k,n} (z^*)^{m}H^kz^n\,,~~
 a_{m,k,n} \in\mC\,.
$$

For technical reasons we consider another type of the generators  $x=Hz, ~~x^*=z^*H$.
They satisfy the commutation relations
\beq{4.9}
Hx=q^\de xH\,, ~~~x^*H=q^\de Hx^*\,, ~~x^*x=q^2xx^*+q^\de(1-q^2)\,.
\eq
\beq{co}
\ddag f(x^*,H,x)\ddag=\sum_{m,k,n=-\infty}^\infty c_{m,k,n}\ti
w(m,k,n)\,,~~~(c_{m,k,n}\in\mC)\,,
\eq
$$
\ti w(m,k,n)=(x^*)^mH^kx^n \,.
$$
In the definition of ${\bf L}_{\de,q}$ we assume that (\ref{co}) are the
formal series. We also consider the "self-conjugate" monomials such as $H^\nu$ and
$(x^*x)^\nu$,  where $\nu>0$.

\bigskip
We define {\sl the non-commutative cone} ${\bf C}_q$ as the
associative *-algebra with the unit and the three generators
$$
(\zeta ^*, \al, \zeta), ~~\al^*=\al\,, ~~(\zeta)^*= \zeta ^*
$$
that satisfy the commutation  relations
\beq{4.11}
\al\zeta=\zeta\al, ~~~\zeta^*\al=\al\zeta^*,
~~~\zeta^*\zeta=q^2\zeta\zeta^*.
\eq
Since $\al$ commutes  with $\ze^*, \ze$ we can define
{\sl the quantum absolute}
$\Xi_q$ as the associative *-algebra generated by $(\ze^*, \ze,1)$
\beq{abs}
\ddag f(\ze^*, \ze)\ddag=\sum_{m,n=-\infty}^\infty b_{m,n}(\ze^*)^m \ze^n\,,~~( b_{m,n}\in\mC)\,.
\eq
We formulate  the conditions on  the coefficients in
next Section.

\bigskip
{\bf 2. Quantum Lorentz groups.}

As in the classical case the NLS are related to the quantum Lorentz group.
We consider here the quantum deformation of the
universal enveloping algebra
${\cal U}_q(\SL)$ $(0<q\leq 1)$ and describe a twisted two parameter
family ${\cal U}^{(r,s)}_q(\SL)$ \cite{OR3}.

We start with a pair of the standard ${\cal U}_q(\SLR)$ Hopf algebra.
The first one is generated by $A,B,C,D$ and the unit with
the relations
\beq{4.1}
AD=DA=1\,,~AB=qBA\,,~BD=qDB\,,
\eq
$$
AC=q^{-1}CA\,,~CD=q^{-1}DC\,,~
[B,C]=\frac{1}{q-q^{-1}}(A^2-D^2)\,.
$$

There is a copy of this algebra ${\cal U}^*_q(\SLR)$ generated by
$A^*,B^*,C^*,D^*$ with the relations following from (\ref{4.1}).
The star
generators commute with $A,B,C,D$. This algebra is the Hopf
algebras with the coproduct
$$
\D(A)=A\ot A\,,
$$
\beq{4.2}
\D(B)=(A^*)^{-r}A\otimes B+B\otimes D(A^*)^s\,,
\eq
$$
\D(C)=(A^*)^rA\otimes C+C\otimes D(A^*)^{-s}\,,
$$
with the counit
$$
\ve\mat{A}{B}{C}{D}=\mat{1}{0}{0}{1}\,,
$$
and the antipode
$$
S\mat{A}{B}{C}{D}=
\mat{D}{-q^{-1}(A^*)^{r-s}B}{-q(A^*)^{s-r}C}{A}\,.
$$

There is the Casimir element in ${\cal U}^{(r,s)}_q(\SL)$ commuting with
any $u\in{\cal U}_q^{r,s}(\SL)$\,.
\beq{4.3}
\Omega_q:=\frac{(q^{-1}+q)(A^2+A^{-2})-4}{2(q^{-1}-q)^2}+
\frac{1}{2}(BC+CB)\,.
\eq
In what follows we put $r=0$.

It was proven in Ref.\,\cite{OR3} that
${\bf L}_{\de,q}$ is a right ${\cal U}^{(0,s)}_q(\SL)$-module.
The right actions of $A, B, C, A^*$ on generators $z^*, H, z$
take the form
$$
w(m,k,n).A=q^{-n+\frac{k}{2}}w(m,k,n)\,,
~~~~~~~w(m,k,n).A^*=q^{\frac{(1-\de)(-2m+k)}s}w(m,k,n)\,,
$$
\beq{4.7}
w(m,k,n).B=q^{-n+\frac{k+1}{2}}
\frac{1-q^{2n}}{1-q^2}w(m,k,n-1)\,,
\eq
$$
w(m,k,n).C=q^{n-\frac{3(k-1)}2+\de(k-1)}
\frac{1-q^{2m}}{1-q^2}w(m-1,k-2,n)-
q^{-n+\frac{k+3}{2}}
\frac{1-q^{2n-2k}}{1-q^2}w(m,k,n+1)\,,
$$
and it follows from (\ref{4.3}) (\ref{4.7}) that
$$
w(m,k,n).\Om_q=q^{-k+1}\frac{(1-q^{k+1})^2}{(1-q^2)^2}w(m,k,n)+
$$
\beq{4.8}
+q^{(\de-1)(k-1)}\frac{(1-q^{2m})(1-q^{2n})}{(1-q^2)^2}w(m-1,k-2,n-1)\,.
\eq
or
$$
f.\Om_q =\frac{ q^{-1}f(q^{-1}x^*,q^{-1}H,q^{-1}x)
+qf(qx^*,qH,qx)}{(1-q^2)^2}
$$
$$
-(1-q^2)^2q^{\de-1}\p_{x^*}\p_xf(q^{-1}x^*,q^{\de-1}H,q^{-1}x)\,.
$$

Note, that when $q \to 1 ~~\Om_q \to \Om$ (\ref{Om}).

It follows from (\ref{4.9}), (\ref{4.7}), and (\ref{4.8}) that the
actions of $A, B, C, A^*$ on $\tilde w(m,k,n)=(x^*)^mH^kx^n$ have the form
$$
\tilde w(m,k,n).A=q^{\frac{m+k-n}{2}}\tilde w(m,k,n)\,,
~~~~~\tilde w(m,k,n).A^*=q^{\frac{(1-\de)(-m+k+n)}s}\tilde w(m,k,n)\,,
$$
$$
\ti w(m,k,n).B=q^{\frac{m+k-n+1}{2}-\de(n-1)}
\frac{1-q^{2n}}{1-q^2}\tilde w(m,k+1,n-1)\,,
$$
$$
\tilde w(m,k,n).C=q^{-\frac{(3m+3k+n-3)}2+\de(k+n)}
\frac{1-q^{2m}}{1-q^2}\tilde w(m-1,k-1,n)+
$$
$$
+q^{-\frac{3m+3k+n-3}{2}+\de n}
\frac{1-q^{2m+2k}}{1-q^2}\tilde w(m,k-1,n+1)\,,
$$
and
$$
\tilde{w}(m,k,n).\Om_q=q^{-m-k-n+1}
\frac{(1-q^{m+k+n+1})^2}{(1-q^2)^2}\tilde{w}(m,k,n)+
$$
\beq{4.10}
+q^{-m-k-n+1+\de(k+1)}
\frac{(1-q^{2m})(1-q^{2n})}{(1-q^2)^2}\tilde{w}(m-1,k,n-1)\,.
\eq

The non-commutative cone ${\bf C}_q$ is also the right module.
We define the actions of $A, B, C, A^*$ on
$\tilde v(m,k,n)=(\zeta^*)^m\al^k\zeta^n$ that compatible with the coproduct
in ${\cal U}_q^{(0,s)}(\SL)$
$$
v(m,k,n).A=q^{\frac{m+k-n}{2}}v(m,k,n),
~~~~~v(m,k,n).A^*=q^{\frac{-m+k+n}s}v(m,k,n)\,,
$$
$$
v(m,k,n).B=q^{\frac{m+k-n+1}2}
\frac{1-q^{2n}}{1-q^2}v(m,k+1,n-1)\,,
$$
$$
v(m,k,n).C=q^{-\frac{3m+3k+n-3}2}
\frac{1-q^{2m+2k}}{1-q^2}v(m,k-1,n+1)\,,
$$
and
$$
v(m,k,n).\Om_q=q^{-m-k-n+1}\frac{(1-q^{m+k+n+1})^2}{(1-q^2)^2}v(m,k,n)+
$$
\beq{4.12}
+q^{-m-k-n+1}\frac{(1-q^{2m})(1-q^{2n})}{(1-q^2)^2}v(m-1,k,n-1)\,.
\eq

It follows from (\ref{4.12}) that $\Om_q$ acts on $v(m,k,n)$ in the
same way as on $w(m,k,n)$ for $\de=0$ and $\al$ playing the role of $H$.

\section{The $q$-Fourier transform and the functional spaces}
\setcounter{equation}{0}

Consider the algebra $\widetilde{\Xi}_q$ generated by $(\xi^*,\xi)$
with the commutation relation
\beq{xi}
\xi\xi^*=q^{2}\xi^*\xi
\eq
and the formal series
\beq{5.2}
\psi(\xi^*,\xi)=\sum_{m=-\infty}^\infty a_{m,n}(\xi^*)^m\xi^n\,.
\eq
Let $L$ be the map of $\widetilde{\Xi}_q$ to the space of functions
on the two-dimensional lattice $C(\mZ\oplus\mZ)$
\beq{f}
L~:~\widetilde{\Xi}_q\to C(\mZ\oplus\mZ)
\eq
\beq{5.3}
L(\psi)=\sum_{m=-\infty}^\infty a_{m,n}q^{2mk}q^{2nl}=
\psi_{k,l}\,,
\eq
and ${\bf K}$ be the algebra of
 the functions (\ref{5.2}) such that the series (\ref{5.3})
converge absolutely for any $k,l\in\mathbb{Z}$.

\begin{defi}\label{d5.1}
The factor space $\widehat{{\bf K}}={\bf K}/\ker  L$ is called the skeleton space
and $L$ (\ref{f}) is the skeleton map.
\end{defi}
In what follows we deal with the skeleton space only.
We can define the skeleton space in the one-dimensional case as
well.

\begin{defi}\label{d5.2}
$\psi\in \widehat{{\bf K}}$ (\ref{5.2}) is a finite function if $\psi_{k,l}=0$
for any $k<-K$ or $l<-L$ and $K, L$ are some positive integers.
\end{defi}
Let $\clK$ be the subalgebra of finite functions from ${\bf K}$.
Define the $q$-Fourier transform on $\clK$
\beq{5.4}
(\clF^{-1}\psi)(\ze^*,\ze)=\phi(\ze^*,\ze)=
\frac1{4\Theta^2_0}\int\int d_{q^2}\xi^*d_{q^2}\xi \psi(\xi^*,\xi)\,,
 \bfE(-q^2\xi^*\ze^*) \bfE(-q^2\xi\ze)\,.
\eq
where the integral is the Jackson integral (\ref{2.7}) and
$\Theta_0$ is determined by (\ref{2.4}).
The $q^2$-integral is well defined because $\psi(\ze^*,\ze) \in\clK $.
The inversion formula has the form
\beq{5.5}
(\clF\phi)(\xi^*,\xi)=\psi(\xi^*,\xi)=\int\int d_{q^2}\ze^*
 \bfe(\xi^*\ze^*)\phi(\ze^*,\ze) \bfe(\xi\ze)  d_{q^2}\ze\,.
\eq
 Let $\clZ$ be the image (\ref{5.4}) of $\clK$. It follows from
the last relation and (\ref{xi}) that $(\ze^*,\ze)$ can be identified
with the absolute generators (\ref{4.11}), and therefore  $\clZ$ is a
subalgebra of $\Xi_q$.

\begin{predl}\label{p5.1} The maps
$$
\clF^{-1}\circ\clF: {\clZ}\rightarrow{\clZ}\,,
$$
$$
\clF\circ\clF^{-1}: {\clK}\rightarrow{\clK}
$$
are the identity maps on ${\clK}$ and ${\clZ}$
correspondingly.
\end{predl}
{\bf Proof}.  For brevity we consider the one-dimensional case.
To prove Proposition we consider the value of the Fourier
transform on the $q^2$-lattice and show that
\beq{5.6}
\clF\circ\clF^{-1}(\psi(q^{2n}))=\psi(q^{2n})\,,
\eq
\beq{5.7}
\clF^{-1}\circ\clF(\phi(q^{2n}))=\phi(q^{2n})\,.
\eq
Consider the first relation
$$
\clF\circ\clF^{-1}(\psi(\xi))=\frac{1-q^2}{2\Theta_0}\int
\bfe(\xi\ze)\int d_{q^2}u\psi(u)\bfE(-q^2u\ze)d_{q^2}\ze=
$$
$$
\frac{1-q^2}{2\Theta_0}q^{2n}\int d_{q^2}u\psi(u)\int
\bfe(\xi\ze)\bfE(-q^2u\ze)d_{q^2}\ze\,.
$$
It follows from Lemma A.1 that we come to (\ref{5.6}).
(\ref{5.7}) is proving just in the same way. \rule{5pt}{5pt}
\bigskip

Now we construct the Fourier transform on ${\mathbf L}_{\de,q}$ with respect to the
"horospheric" generators $(x^*,x)$ in the similar way as above.
Let $f(x^*,H,x)$ be an element from ${\mathbf L}_{\de,q}$ (\ref{co}),
such that the inverse Fourier integral
\beq{5.8}
\clF^{-1}(f)(y^*,H,y)=
\frac1{4\Theta_0^2}\int\int
d_{q^2}x^*\ddag f(x^*,H,x)\bfE(-q^2y^* x^*)\bfE(-q^2yx)\ddag d_{q^2}x
\eq
is well defined.
We preserve the notion ${\mathbf L}_{\de,q}$ for the space of these functions
and define the algebra
\beq{tl}
\widetilde{\mathbf L}_{\de,q}=\clF^{-1}({\mathbf L}_{\de,q})
\eq
with the generators $(y^*,H,y)$ and the commutation relations
\beq{cre}
yH=Hy\,,~~~y^*H=Hy^*\,,
\eq
$$
yy^*=q^{-2}y^*y[1+q^\de(q^2-1)y^*y]^{-1}\,.
$$
The direct Fourier transform takes the form
\beq{5.9}
(\clF g)(x^*,H,x)=\int\int\bfe(y^* x^*)
d_{q^2}y^*\ddag g(y^*,H,y)\ddag d_{q^2}y\bfe(yx)\,.
\eq
Then as before
$$
\clF \circ\clF^{-1}=Id ~\hbox{on}~{\mathbf L}_{\de,q}\,,~~~
\clF^{-1}\circ\clF=Id ~\hbox{on}~\widetilde{\mathbf L}_{\de,q}\,.
$$

Let $\nu\geq 0$ and $\widetilde{\clW}_\nu$ be the space of functions with
a fixed singularity.
It is constructed by means of a pair elements $g_1,g_2\in\widetilde{\mathbf L}_{\de,q}$ as follows
\beq{cW}
\widetilde{\clW}_\nu =\left\{\begin{array}{rcl}
g_1(y^*,H,y)+(y^*)^\nu g_2(y^*,H,y)y^\nu &{\rm if} & \nu\ne n\in\mathbb{N}\,,\\
g_1(y^*,H,y)+\ln y^*g_2(y^*,H,y)+g_2(y^*,H,y)\ln y &{\rm if} & \nu=n\in\mathbb{N}\,.\\
\end{array}\right\}
\eq
We define $\clW_\nu$ as the image of $\widetilde{\clW}_\nu$
   by the Fourier transform (\ref{5.9})
\beq{cw}
\clW_\nu = \clF(\widetilde{\clW}_\nu )\,.
\eq

\section{The Poisson kernel}
\setcounter{equation}{0}

The Poisson kernel is the element of the algebra
${\bf L}_{q,\de}\otimes\Xi_q$ determined by the series
\beq{7.1}
\clP_\nu((x^*\otimes1-1\otimes\ze^*)\,,H\,,
(1\otimes x-\ze\otimes1))=
\eq
$$
=\sum_{k=0}^\infty(-1)^k\frac{(q^{2\nu+2},q^2)_k}{(q^2,q^2)_k}
q^{(2-\nu\de-2\de)k}(x^*\otimes1-1\otimes\ze^*)^kH^{\nu+1}
(-\ze\otimes1+1\otimes x)^k\,.
$$
Let
$$
\Om_q^\nu=\Om_q-q^{-\nu+2}\left(\frac{1-q^\nu}{1-q^2}\right)Id\,.
$$
and $ F_\nu(x^*,H,x)\in {\bf L}_{q,\de}$ be a solution to the equation
\beq{7.6}
F_\nu(x^*,H,x).\Om^\nu_q=0 \,.
\eq
We can formulate now our main result.
\begin{predl}\label{p8.2}
For any $\phi\in\clZ $ on the absolute $\Xi_q$\\
$\bullet$ The function
\beq{8.2}
F_\nu(x^*,H,x)=(\clP_\nu*\phi)(x^*,H,x)\,,
\eq
is a solution to the equation (\ref{7.6});\\
$\bullet \bullet$ $F_\nu(x^*,H,x)\in\clW_\nu$ (\ref{cw}).
\end{predl}
Here the convolution is defined as
\beq{8.1}
(\clP_\nu*\phi)(x^*,H,x)=
\int\int d_{q^2}\ze^*\ddag
\clP_\nu((x^*\otimes1-1\otimes\ze^*)\,,H\,,
(-\ze\otimes1+1\otimes x))\phi(\ze^*,\ze)
\ddag d_{q^2}\ze\,.
\eq
We postpone the proof of this statement to last Section and formulate
here some intermediate steps.

Let $\si$ be the involution $\si f(x^*,H,x)= f(-x^*,H,-x)$.
\begin{predl}\label{p7.1}
The Poisson kernel has the integral representation
\beq{7.2}
\clP_\nu(x^*,H,x)=\si\clF(Q_\nu)\,,
\eq
where
\beq{7.3}
Q_\nu(y^*,H,y)=
\eq
$$
=\frac{1+q}{2\Te_0^2}\frac{q^{\nu^2+\nu+\de(\oh\nu^2+\nu)}}{\G_{q^2}(\nu+1)}
(y^*)^{\frac\nu2}H^{\frac\nu2}
\ddag K_\nu^{(2)}(2(y^*)^\oh y^\oh q^\de(1-q^2);q^2)H\ddag
H^{\frac\nu2}y^{\frac\nu2}\,,
$$
$K_\nu^{(2)}$ is the $q^2$-Bessel-Macdonald function of kind $2$
\cite{R,OR1}, and $\G_{q^2}(\nu+1)$ is the\\ $q^2-\G$-function (\ref{G}).
\end{predl}

Consider the Fourier transform of the left hand side (\ref{7.6})
$$
\si\clF^{-1}\left( f.\Om^\nu_q\right)=
\si\clF^{-1}(f).\widetilde{\Om}^\nu_q\,,~~~\widetilde{\Om}^\nu_q =\clF^{-1}\Om^\nu_q \,.
$$
For $ \si\clF^{-1}(f)= g(y^*,H,y)$ we have the equation
$ g(y^*,H,y).\widetilde{\Om}^\nu_q=0$, or
$$
q^{-1}g(q^{-1}y^*,qH,q^{-1}y)-(q^\nu+q^{-\nu})g(y^*,H,y)+
qg(qy^*,q^{-1}H,qy)=
$$
\beq{7.10}
=(1-q^2)^2q^{\de+1}y^*g(qy^*,q^{\de-1}H,qy)y\,.
\eq

The $q^2$-Fourier transform of the $q^2$-Poisson kernel satisfies (\ref{7.10})
\beq{Q}
Q_\nu (y^*,H,y).\widetilde{\Om}^\nu_q=0\,.
\eq
The statement is verified directly using the series
representation of the $q^2$-Bessel-Macdonald function.

\begin{predl}\label{p7.4}
If $g(y^*,H,y)$ is a solution of (\ref{7.10}), then the
product\\
$g(y^*,H,y)\psi(\al y^*,\al y)$ is a solution to the same equation for
any function $\psi(\al y^*,\al y)$.
\end{predl}

Let
\beq{QQ}
{\cal Q}_\nu(y^*,H,y;\xi^*,\xi)=Q_\nu(y^*,H,y )
\de(y^*\otimes 1-1\otimes\xi^*, y\otimes 1-1\otimes\xi)\,.
\eq
Here delta-function is the kernel of the integral transform
$\widetilde{\Xi}_q\to\widetilde{\mathbf L}_{\de,q}$
$$
h(y^*,y)=\int d_{q^2}\xi^*\ddag h(\xi^*,\xi)
\de(y^*\otimes 1-1\otimes\xi^*, y\otimes 1-1\otimes\xi)
\ddag d_{q^2}\xi\,.
$$
Thereby, the multiplication by $Q_\nu$ carries out this map.

\begin{predl}\label{p7.5}
The operator of multiplication on $Q_\nu(y^*,H,y)$ transforms the space $\clK$
into the space $\widetilde{\clW}_\nu$.
\end{predl}

We illustrate our construction by the following commutative diagram
$$
\def\normalbaselines{\baselineskip20pt
     \lineskip3pt   \lineskiplimit3pt}
\def\mapright#1{\smash{
     \mathop{\longrightarrow}\limits^{#1}}}
\def\mapleft#1{\smash{
     \mathop{\longleftarrow}\limits^{#1}}}
\def\mapdown#1{\Big\downarrow\rlap
{$\vcenter{\hbox{$\scriptstyle#1$}}$}}
\matrix{(\psi\in\clK,~\widetilde{\Xi}_q)  &\mapright{\clF^{-1}}  & (\phi\in\clZ ,~\Xi_q) \cr
\mapdown{Q_\nu\times}&                 &\mapdown{\clP_\nu *}\cr
(\widetilde{\clW}_\nu,~ \widetilde{\mathbf L}_{\de,q})&
\mapright{\si\clF}                     &
(F_\nu\in\clW_\nu,~ {\mathbf L}_{\de,q}) \cr}
$$
It means that one can start with a finite function $\psi\in\clK$ and then come to the
solution $F_\nu$ by one of the two possible ways.

\section{The $q^2$-Fourier transform of spherical symmetric functions}
\setcounter{equation}{0}

\begin{defi}\label{d6.1} A function $f\in{\mathbf L}_{\de,q}$ is
{\sl spherical symmetric} if it depends on the product of $x^*x$
\beq{6.1}
f(x^*,H,x)=\sum_{l,k=-\infty}^\infty c_{l,k}(x^*)^lH^kx^l\,.
\eq
\end{defi}

Any ordered element from ${\mathbf L}_{\de,q}$ can be represented
in the form \cite{SSV}
\beq{6.2}
f(x^*,H,x)=\sum_{r=1}^\infty (x^*)^r\phi_{-r}(x^*,H,x)+
\phi_0(x^*,H,x)+\sum_{r=1}^\infty\phi_r(x^*,H,x)x^r\,,
\eq
where $\phi_r(x^*,H,x)$ are spherically symmetric.

Consider the inverse Fourier transform of the  spherical
symmetric function (\ref{6.1})
$$
\clF^{-1}f=g(y^*,H,y)=\frac1{4\Te_0^2}\int\int d_{q^2}x^*
\ddag\bfE(q^2y^* x^*)f(x^*,H,x)\bfE(q^2yx)\ddag d_{q^2}x\,.
$$
Using (\ref{2.2})  we obtain
$$
g(y^*,H,y)=
$$
$$
=\frac1{4\Te_0^2}(1-q^2)\sum_{n=-\infty}^\infty q^{2n}
\sum_{m=0}^\infty\frac{(-1)^m(1-q^2)^{2m}q^{2m(m+1)}q^{2nm}}{(q^2,q^2)_m^2}
(y^*)^m\sum_{l,k}c_{lk}q^{2nl}H^ky^m=
$$
$$
=\frac1{4\Te_0^2}(1-q^2)\sum_{n=-\infty}^\infty q^{2n}
\ddag J_0^{(2)}(2(y^*)^\oh y^\oh q^nq(1-q^2);q^2)
\sum_{l,k}c_{lk}q^{2nl}H^k\ddag\,,
$$
where $J_0^{(2)}$ is the $q^2$-Bessel function of kind 2 \cite{GR}.
Using (\ref{2.2a}) the last expression can be rewritten as the integral
\beq{6.4}
(\clF^{-1}f)(y^*,H,y)=\frac{1+q}{4\Te_0^2}\int_0^\infty
\ddag J_0^{(2)}(2(y^*)^\oh y^\oh q(1-q^2)\rho;q^2)
\sum_{l,k}c_{lk}\rho^{2l}H^{k-2}\ddag\rho d_q\rho
\eq
$$
=\frac{1+q}{4\Te_0^2}\int_0^\infty
\ddag J_0^{(2)}(2(y^*)^\oh y^\oh q(1-q^2)\rho;q^2) f(x^*,H,x)|_{x^*=\rho,x=\rho}
\ddag\rho d_q\rho\,.
$$

The inversion formula has the form
\beq{6.5}
\clF (g)=f(x^*,H,x)=
\eq
$$
=(1+q)\int_0^\infty\ddag J_0^{(1)}(2(x^*)^\oh x^\oh(1-q^2)r;q^2)
\sum_{l,k}b_{lk}r^{2l}H^k\ddag rd_qr
$$
$$
=(1+q)\int_0^\infty\ddag J_0^{(1)}(2(x^*)^\oh x^\oh(1-q^2)r;q^2)
g(y^*,H,y)|_{y^*=r,y=r} \ddag rd_qr \,.
$$
Here $J_0^{(1)}$ is the $q^2$-Bessel function of kind 1 \cite{GR}.

\section{Proofs}
\setcounter{equation}{0}

{\bf 1. Proof of \ref{p7.1}}.

Consider the power series
\beq{7.4}
\clP_\nu(x^*,H,x)=\sum_{k=0}^\infty(-1)^k
\frac{(q^{2\nu+2},q^2)_k}{(q^2,q^2)_k}
q^{(2-\nu\de-2\de)k}(x^*)^kH^{\nu+1}x^k\,.
\eq
and find its Fourier transform $Q_\nu(y^*,H,y)$ as the
function from $\widetilde{\mathbf L}_{\de,q}$.
Since $\clP_\nu(x^*,H,x)$   is the spherically symmetric we have,
following (\ref{6.4})
$$
Q_\nu(y^*,H,y)=\frac{1+q}{4\Te_0^2}
\int_0^\infty\ddag J_0^{(2)}(2(y^*)^\oh y^\oh q(1-q^2)\rho;q^2)
\sum_ka_kq^{(2-\nu\de-2\de)k}\rho^{2k}H^{\nu+1}\ddag\rho d_q\rho\,,
$$
where $a_k=\frac{(-1)^k(q^{2\nu_2},q^2)_k}{(q^2,q^2)_k}$.
On the other hand one can find that
\beq{7.5}
\sum_{k=0}^\infty
a_kq^{(2-\nu\de-2\de)k}\rho^{2k}=\frac{(-q^{(2+\nu)(2-\de)}\rho^2,q^2)_\infty}
{(-q^{2-(2+\nu)\de}\rho^2,q^2)_\infty}\,.
\eq
Hence,
$$
Q_\nu(y^*,H,y)=\frac{1+q}{4\Te_0^2}\times
$$
$$
\times\int_0^\infty\rho\frac{(-q^{(2+\nu)(2-\de)}\rho^2,q^2)_\infty}
{(-q^{2-(2+\nu)\de}\rho^2,q^2)_\infty}\ddag H^{\nu+1}
J_0^{(2)}(2(y^*)^\oh y^\oh q(1-q^2)\rho;q^2)\ddag d\rho\,.
$$

It follows from Ref.\,\cite{OR2,R} that the last expression is the integral
representation of the $q^2$-Bessel-Macdonald function, i.e.
\beq{7.6a}
Q_\nu(y^*,H,y)=B(y^*)^{\frac\nu2}H^{\frac\nu2}
\ddag K_\nu^{(2)}(2(y^*)^\oh y^\oh q^\de(1-q^2);q^2)H\ddag
H^{\frac\nu2}y^{\frac\nu2}\,.
\eq
To calculate $B$ we restrict the integral to the common
kernel  $y^*=0,~y=0$. Then we obtain \cite{R,OR1}
$$
\frac{1+q}{4\Theta_0^2}\int d_q\rho^2\frac{(-q^{(2+\nu)(2-\de)}\rho^2,q^2)_\infty}
{(-q^{2-(2+\nu)\de}\rho^2,q^2)_\infty}=
B\oh q^{-\nu^2+\nu-\de(\oh\nu^2+\nu)}\G_{q^2}(\nu)\,.
$$

Now we calculate the integral in the right side. Note that
$$
\p_{\rho^2}E_{q^2}(a\rho^2)=\frac a{1-q^2}E_{q^2}(aq^2\rho^2),
~~~\p_{\rho^2}e_{q^2}(b\rho^2)=\frac b{1-q^2}e_{q^2}(b\rho^2)\,.
$$
Hence
$$
\int_0^\infty d_q\rho^2e_{q^2}(-q^{2-(2+\nu)\de}\rho^2)
E_{q^2}(q^{(2+\nu)(2-\de)}\rho^2)=
$$
$$
=(1-q^2)q^{2-(2+\nu)(2-\de)}\int_0^\infty d_q\rho^2\p_{\rho^2}
E_{q^2}(q^{(2+\nu)(2-\de)-2}\rho^2)e_{q^2}(-q^{2-(2+\nu)\de}\rho^2)=
$$
$$
=(1-q^2)[\lim_{m\to\infty}\frac{-q^{(2+\nu)(2-\de)-2m},q^2)_\infty}
{-q^{2-(2+\nu)\de-2m},q^2)_\infty}-1]-
$$
$$
=-(1-q^2)q^{2-(2+\nu)(2-\de)}\int_0^\infty d_q\rho^2
E_{q^2}(q^{(2+\nu)(2-\de)}\rho^2)
\p_{\rho^2}e_{q^2}(-q^{2-(2+\nu)\de}\rho^2)=
$$
$$
=-(1-q^2)+q^{-2\nu}\int_0^\infty d_q\rho^2\p_{\rho^2}
E_{q^2}(q^{(2+\nu)(2-\de)}\rho^2)e_{q^2}(-q^{2-(2+\nu)\de}\rho^2)\,.
$$
Then we have
$$
\int_0^\infty d_q\rho^2e_{q^2}(-q^{2-(2+\nu)\de}\rho^2)
E_{q^2}(q^{(2+\nu)(2-\de)}\rho^2)=-\frac{1-q^2}{1-q^{-2\nu}}\,.
$$
It implies that
\beq{7.9}
B=\frac{1+q}{2\Theta_0^2\G_{q^2}(\nu+1)}q^{\nu^2+\nu+\de(\oh\nu^2+\nu)}
\eq
and the representation (\ref{7.3}).\rule{5pt}{5pt}

\bigskip
{\bf 2. Proof of 4.3}

Consider the equation (\ref{7.10}) for the product
$g(y^*,H,y)\psi(\al y^*,\al y)$.
Note that $\al$ plays the role of $H$ and the action of the generators
is the same as for $\de=0$ (compare (\ref{4.10}) and (\ref{4.12})).
In this way we have
$$
q^{-1}g(q^{-1}y^*,qH,q^{-1}y)\psi(q\al q^{-1}y^*,q\al q^{-1}y)-
(q^\nu+q^{-\nu})g(y^*,H,y)\psi(\al y^*,\al y)+
$$
$$
qg(qy^*,q^{-1}H,qy)\psi(q^{-1}\al qy^*,q^{-1}\al qy)=
$$
$$
(1-q^2)^2q^{\de+1}y^*g(qy^*,q^{\de-1}H,qy)\psi(q^{-1}\al qy^*,q^{-1}\al qy)y\,.
$$
\rule{5pt}{5pt}

\bigskip
{\bf 3. Proof of 4.4}

We can extract the series expansion for $K_\nu^{(2)}$ from the series expansion
of $I_\nu^{(2)}$ and $I_{-\nu}^{(2)}$ for $\nu\ne n\in\mathbb{N}$ \cite{OR1}.
It gives us the following expansion
$$
Q_\nu=\frac{(1+q)(1-q^2)}{4\Theta^2(1-q^{2\nu})}q^{2\nu}\times
$$
$$
\times[\sum_{l=0}^\infty\frac{q^{2l(l-\nu)}(1-q^2)^{2l}}
{(q^2,q^2)_l(q^{-2\nu+2},q^2)_l}
q^{\de l(\nu+2)}(y^*)^lH^{\nu+1}y^l-
$$
$$
-\frac{\G_{q^2}(1-\nu)}{\G_{q^2}(1+\nu)}
\sum_{l=0}^\infty\frac{q^{2l(l+\nu)}(1-q^2)^{2l}}
{(q^2,q^2)_l(q^{2\nu+2},q^2)_l}
q^{\de l(\nu+2)}(y^*)^{\nu+l}H^{\nu+1}y^{\nu+l}]=
$$
$$
=\Psi_1+(y^*)^\nu\Psi_2 y^\nu.
$$
The functions $\Psi_1$ and $\Psi_2$ are multipliers in $\clK$ and we
come to the statement of Proposition for $\nu\neq\mN$.

Using the expression of $K_n^{(2)}$ \cite{OR1} it can be proved similarly
that $Q_n $ has the logarithmic singularity (\ref{cW}). \rule{5pt}{5pt}

\bigskip
{\bf 4. Proof of 4.1}

Let $\clF\phi=\psi$. Remind that $\phi$ and $\psi$ depend on
$\al$ (see n.4). Taking into account (\ref{7.1})
$$
\clP_\nu=\clF[\bfe(\ze^*y^*)Q_\nu(y^*,H,y)\bfe(y\ze)]\,,
$$
we have
$$
\si\clF^{-1}(\clP_\nu*\phi)=\si\clF^{-1}(\clF[\bfe(\ze^*y^*)
Q_\nu(y^*,H,y)\bfe(y\ze)]*\clF^{-1}\psi)=
$$
$$
\frac1{4\Theta^2}\int\int d_{q^2}u^*Q_\nu(u^*,H,u)d_{q^2}u
\int\int d_{q^2}v^*\psi(v^*,\al,v)d_{q^2}v
\int\int d_{q^2}x^*\bfE(q^2x^*y^*)\bfe(-x^*u^*)\times
$$
$$
\times\bfe(-ux)\bfE(q^2yx)d_{q^2}x\times
$$
$$
\times\int\int d_{q^2}\ze^*\bfe((1\otimes\ze^*)u^*)
\bfE(-q^2v^*(\ze^*\otimes1))\bfE(-q^2(1\otimes\ze)v)
\bfe(u(\ze\otimes1))d_{q^2}\ze\,.
$$

It follows from Lemma A.1 that
$$
\int\bfe(-ux)\bfE(q^2yx)d_{q^2}x=
\left\{\begin{array}{lcr}
\frac{2\Theta_0}{1-q^2}& {\rm for}&u=y\\
0 & {\rm for} & u\ne y
\end{array}
\right.
$$
and
$$
\int\bfE(-q^2(1\otimes\xi)v)
\bfe(u(\ze\otimes1)d_{q^2}\ze=
$$
$$
=\left\{\begin{array}{lcr}
\frac{2\Theta_0)}{1-q^2}&
{\rm for}&u\otimes1=1\otimes v \\
0 & {\rm for} & u\otimes1\ne1\otimes v
\end{array}
\right.
$$
Hence
$$
\si\clF(\clP_\nu*\phi)=Q_\nu(y^*,H,y)\psi(y^*,\al,y)\,.
$$
\rule{5pt}{5pt}

\section*{Appendix A. $q$-relations}
\setcounter{equation}{0}
\def\theequation{A.\arabic{equation}}

We assume that $|q|<1$. Let us recall some notations \cite{GR}.
We consider the $q^2$-exponentials
\beq{2.1}
e_{q^2}(z)=\sum_{n=0}^\infty\frac{z^n}{(q^2,q^2)_n}=
\frac1{(z,q^2)_\infty},\qquad|z|<1,
\eq
\beq{2.2}
E_{q^2}(z)=\sum_{n=0}^\infty\frac{q^{n(n-1)}z^n}{(q^2,q^2)_n}=
(-z,q^2)_\infty,
\eq
and $q^2-\G$-function
\beq{G}
\G_{q^2}(\nu)=\frac{(q^2,q^2)_\infty}{(q^{2\nu},q^2)_\infty}(1-q^2)^{1-\nu}.
\eq
Introduce the notion of the $q^2$-exponentials defined on a tensor product
\beq{2.3}
\bfe(y\ze)=e_{q^2}(i(1-q^2)(y\otimes \ze)),
~~~~\bfE(y\ze)=E_{q^2}(i(1-q^2)(y\otimes \ze)).
\eq
Consider
$$
{\bf Q}(z,q)=(1-q^2)\sum_{m=-\infty}^\infty\frac1
{zq^{2m}+z^{-1}q^{-2m}}
$$
and let
\beq{2.4}
\Theta_0={\bf Q}(1-q^2,q).
\eq

Let ${\cal A}=C[z,z^{-1}]$ be the algebra of formal Laurent series.
The $q^2$-derivative of a function $f(z)\in {\cal A}$ is defined as follows
\beq{2.5}
\p_zf(z)=(f(z)-f(q^2z))\frac{z^{-1}}{1-q^2}.
\eq

The functions (\ref{2.3}) satisfy conditions:
$$
\p_\xi\bfE(y\ze)=(y\otimes1)\bfE(q^2y\ze),
$$
$$
\p_y\bfE(y\ze)=(1\otimes\ze)\bfE(q^2y\ze),
$$
$$
\p_\ze\bfe(y\ze)=(y\otimes1)\bfe(y\ze),
$$
\beq{2.6}
\p_y\bfe(y\ze)=(1\otimes\ze)\bfe(y\ze).
\eq

The $q^2$-integral (Jackson integral \cite{GR}) is defined to be the
following map $I_{q^2}$ of the algebra ${\cal A}$ into the space of formal
numerical series:
\beq{2.7}
I_{q^2}f=\int d_{q^2}zf(z)=(1-q^2)\sum_{m=-\infty}^\infty q^{2m}
[f(q^{2m})+f(-q^{2m})].
\eq
It follows from this definition that
\beq{2.8}
I_{q^2}\p_zf(z)=0.
\eq
If the series in right side of (\ref{2.7}) is nonconvergent
(\ref{2.8}) is a regularization of nonconvergent $q^2$-integral.

We need also another type of the Jackson integral
\beq{2.2a}
\int_0^\infty d_qxf(x)=(1-q^2)\sum_{m=-\infty}^\infty q^{m}f(q^{m})\,.
\eq

The function $f(z)$ is absolutely $q^2$-integrable if the series
$$
\sum_{m=-\infty}^\infty q^{2m}[|f(q^{2m})|+|f(-q^{2m})|]
$$
converges.

{\bf Lemma A.1}
\beq{2.9}
\int\bfE(q^2y\ze)\bfe(-u\ze)d_{q^2}\ze=
\left\{\begin{array}{lcl}
\frac2{1-q^2}\Theta_0u^{-1} &{\rm if}& ~~~y=u\\ 0 &{\rm if}& ~~~y\ne u\,,\\
\end{array}
\right.
\eq
\beq{2.10}
\int d_{q^2}y\bfE(q^2y\xi)\bfe(-y\ze)=
\left\{\begin{array}{lcl}
\frac2{1-q^2}\Theta_0\xi^{-1} & {\rm if} & ~~\ze=\xi\\
0 & {\rm if} & ~~\ze\ne \xi\,,\\
\end{array}
\right.
\eq
{\sl where} $\Theta_0$ {\sl determined by (\ref{2.4}).}

{\bf Proof.}
It follows from (\ref{2.6}) that
$$
\p_\ze[\bfE(y\ze)\bfe(-u\ze)]=i((y-u)\otimes1)\bfE(q^2y\ze)\bfe(-u\ze)\,.
$$
Hence, if $y\ne u$, then
$$
(y-u)\int\bfE(q^2y\ze)\bfe(-u\ze)d_{q^2}\ze=
\int\p_\ze[\bfE(y\ze)\bfe(-u\ze)]d_{q^2}\ze=0
$$
in accordance with the definition of the $q^2$-integral.

If $y=u$, then
$$
\int\bfE(q^2u\ze)\bfe(-u\ze)d_{q^2}\ze=
u^{-1}\int\bfE(q^2\ze)\bfe(-\ze)d_{q^2}\ze=
$$
$$
=u^{-1}(1-q^2)\sum_{m=-\infty}^\infty q^{2m}
\left[\frac{(-i(1-q^2)q^{2m+2},q^2)_\infty}
{(-i(1-q^2)q^{2m},q^2)_\infty}+
\frac{(i(1-q^2)q^{2m+2},q^2)_\infty}
{(i(1-q^2)q^{2m},q^2)_\infty}\right]=
$$
$$
=u^{-1}(1-q^2)\sum_{m=-\infty}^\infty q^{2m}
\left(\frac1{1+i(1-q^2)q^{2m}}+\frac1{1-i(1-q^2)q^{2m}}\right)=
$$
$$
=u^{-1}2\sum_{m=\infty}^\infty\frac1{(1-q^2)^{-1}q^{-2m}+(1-q^2)q^{2m}}=
u^{-1}\frac2{1-q^2}\Theta_0
$$
(see (\ref{2.4})).

(\ref{2.10}) is proving in just the same way. \rule{5pt}{5pt}

\bigskip
\small{

}

\end{document}